\documentclass[11pt]{amsart}

\usepackage{epsf,amssymb}

\newtheorem{thm}{Theorem}[section]
\newtheorem{lem}[thm]{Lemma}

\newtheorem{prop}[thm]{Proposition}
\newtheorem{rem}[thm]{Remark}

\newcommand{\bbr}{\begin{rem}\em} 
\newcommand{\eer}{\end{rem}}
\newcommand{\bdry}{\partial}

\newcommand{\s}{\vskip.1in}
\def\Z{\hbox{$\mathbb Z$} }

\def\R{\hbox{$\mathbb R$} }

\def\co{\colon\thinspace}

\def\dfn#1{{\em #1}}

\begin{document}

%%%%%%%%%%%%%%%%%%%%%%%%%%%%%%%%%%%%%%%%%%%%%%%%%%%%%%%%%%
\title{On Symplectic Cobordisms}
\date{February 13, 2001}

\author{John B. Etnyre}
\address{Stanford University, Stanford, CA 94305}
\email{etnyre@math.stanford.edu}
\urladdr{http://math.stanford.edu/\char126 etnyre}

\author{Ko Honda}
\address{University of Georgia, Athens, GA 30602\newline
Current Address: IHES, Bures-sur-Yvette, France}
\email{honda@math.uga.edu}
\urladdr{http://www.math.uga.edu/\char126 honda}

\keywords{tight, overtwisted, cobordism, Stein}
\subjclass{Primary 53C15; Secondary 57M50}
\thanks{JE supported by NSF grant DMS-9705949.  KH supported by NSF grant DMS-0072853.}

\begin{abstract}
In this note we make several observations concerning symplectic cobordisms. Among other things we 
show that every contact 3-manifold has infinitely many concave symplectic fillings and that
all overtwisted contact 3-manifolds are ``symplectic cobordism equivalent.''	
\end{abstract}

\maketitle

%%%%%%%%%%%%%%%%%%%%%%%%%%%%%%%%%%%%%%%%%%%%%%%%%%%%%%%%%%

% ************************************************************************
\section{Introduction}
% ************************************************************************

In this note we make several observations concerning (directed) symplectic cobordisms,
Stein cobordisms, and {\em concave} symplectic fillings for contact 3-manifolds. 
Symplectic and Stein cobordisms have recently come to the foreground of symplectic and contact
geometry, largely due to the introduction of a symplectic field theory (SFT) by Eliashberg, 
Hofer and Givental \cite{EGH}.  The goal of SFT is to associate an 
algebraic structure to a given symplectic cobordism.  Though clearly a central notion in 
symplectic and contact geometry, there is surprisingly little concerning symplectic cobordisms 
in the literature. 

We will assume our 3-manifolds are closed and oriented, and our contact 
structures are  oriented and positive.  A contact 3-manifold $(M_1,\xi_1)$ is {\it 
symplectically (resp.\ Stein) cobordant} to another contact manifold $(M_2,\xi_2)$, if there 
exists a symplectic (resp.\ Stein) 4-manifold $(X,\omega)$ with $\bdry X=M_2-M_1$ and a vector 
field $v$ defined on a neighborhood of $(M_1\cup M_2)\subset X$ for which 
$\mathcal{L}_v\omega=\omega$, $v\pitchfork  (M_1\cup M_2)$,
and the normal orientation of $M_1\cup M_2$ agrees with $v$.
We say $(M_1,\xi_1)$ is the {\it concave end} of the cobordism, while 
$(M_2,\xi_2)$ is the {\it convex end}.  We denote the existence of such a 
cobordism by $(M_1,\xi_1)\prec (M_2,\xi_2)$ --- in the 
paper we implicitly assume that $\prec$ refers to a Stein cobordism, unless specified 
otherwise.  Note that symplectic (and Stein) cobordism is 
not an equivalence relation.  For example, a {\it Stein fillable} contact 
structure $(M,\xi)$ (= one satisfying
$\emptyset \prec (M,\xi)$)  cannot be symplectically cobordant to an overtwisted 
contact structure, but the opposite is possible. Our first result is:

\begin{thm}\label{thm1}
Let $(M_1,\xi_1)$ be a contact 3-manifold.  Then there exists a Stein fillable contact
3-manifold $(M_2,\xi_2)$ and a Stein cobordism $(M_1,\xi_1)\prec (M_2,\xi_2)$.
\end{thm}

Though this result indicates the overall structure of the partial order on contact 
3-manifolds induced by cobordisms, there is very little control over the target contact 
manifold $(M_2,\xi_2)$. On the other hand, when $(M_1,\xi_1)$ is overtwisted, there is 
complete freedom in choosing $(M_2,\xi_2)$: 

\begin{thm}   \label{thm2}
Let $(M_1,\xi_1)$ be an overtwisted contact 3-manifold and $(M_2,$ $\xi_2)$ any contact 
3-manifold, tight or overtwisted.  Then there exists a Stein cobordism  
$(M_1,\xi_1)\prec (M_2,\xi_2)$.
\end{thm}

In particular, all overtwisted contact structures are 
equivalent under symplectic or Stein cobordism! 

It is interesting to compare the previous two theorems with recent work of  
Epstein-Henkin \cite{EpsteinH} and de Oliveira \cite{deOliveira} which deal 
with cobordisms between CR-structures. (Here ``CR-structure'' will mean 
``strictly pseudoconvex CR-structure''.)  On any 3-manifold $M$, there is a 1-1 
correspondence between CR-structures and pairs $(\xi,J)$ consisting of a contact 
structure $\xi$ and an almost complex structure $J$ on $\xi$. 
We say a  CR-structure $(\xi, J)$ on $M$ is \textit{fillable}, if there is a 
compact, connected, complex manifold $X$ with $\bdry X=M$, so that the complex 
tangencies to $M$ are $\xi$ and the induced complex structure on $\xi$ is $J.$ 
In \cite{EpsteinH} it was shown that if a CR-manifold $(M_1,\xi_1, J_1)$ is 
Stein cobordant to a fillable CR-manifold $(M_2,\xi_2,J_2)$, then 
$(M_1,\xi_1,J_1)$ is also fillable.  Here we assume Stein cobordisms of 
CR-manifolds respect complex structures. 
Thus, if $(M_1,\xi_1,J_1)\prec (M_2,\xi_2,J_2)$ is a Stein cobordism but $(M_1,\xi_1)$ is not 
Stein fillable, then $(M_2,\xi_2,J_2)$ 
cannot be a fillable CR-structure, even if $(M_2,\xi_2)$ is a Stein fillable 
contact structure. De Oliveira \cite{deOliveira} gave some interesting examples 
of complex (but not Stein) cobordisms from non-fillable CR-structures to 
fillable ones, thus showing the necessity of having a Stein cobordism in the 
Epstein-Henkin result.

Our last result is: 
 
\begin{thm} \label{thm3}
Any contact 3-manifold has infinitely many concave symplectic fillings which are mutually 
non-isomorphic and are not related to each other by a sequence of blow-ups and blow-downs. 
\end{thm}

A {\it convex (resp.\ concave) symplectic filling} of $(M,\xi)$ is a symplectic cobordism 
$(X,\omega)$ from $\emptyset$ to $(M,\xi)$ (resp.\ from $(M,\xi)$ to $\emptyset$).  
The phrase ``symplectic filling,'' without modifiers, is usually reserved for ``convex 
symplectic filling.'' 
Having a (convex) filling is quite restrictive for a contact 3-manifold --- for instance, it 
implies the contact structure is tight.  (Note, however, that there are many tight 
contact structures without such fillings due to Eliashberg \cite{Eliashberg96}, Ding-Geiges 
\cite{DG}, and Etnyre-Honda \cite{EH}.) We 
show that, on the contrary, concave fillings are not restrictive at all. Though this was 
believed for a long time, and specific isolated contact manifolds with infinitely  many 
such fillings are easy to come by,  the degree to which concave fillings are not restrictive is 
perhaps a little surprising.

We assume the reader is more or less familiar with contact geometry and hence we do not include 
any background material here. We refer the reader to \cite{Aebischer} for the basics of contact
geometry, \cite{Eliashberg89} for Lutz twisting, and \cite{EGH, Eliashberg90} for the notions 
of Stein and symplectic cobordisms.

% ************************************************************************
\section{Legendrian surgeries}   \label{sect2}
% ************************************************************************

In this section we give a description of {\it Legendrian surgery}, both on the 3-manifold
level and as a source of Stein filling on the 4-manifold level.  There is some related material 
in \cite{H2} for Legendrian surgeries.

Let $(M,\xi)$ be a contact manifold and $L\subset M$ a closed Legendrian curve.  Let 
$N(L)$ be a {\it standard tubular neighborhood} of the Legendrian curve $L$, with convex 
boundary and two parallel dividing curves.  Choose a framing for $L$ (and a concomitant 
identification $\bdry N(L)\simeq \R^2/\Z^2$) so that the meridian has slope $0$ and the 
dividing curves have slope $\infty$.  With respect to this choice of framing, a Legendrian 
surgery is a $-1$ surgery, where a copy of $N(L)$ is glued to $M\setminus N(L)$ so that the new 
meridian has slope $-1$.  Here, even though the boundary characteristic foliations may not 
exactly match up a priori,  we use Giroux's Flexibility Theorem 
\cite{Giroux91, H1} and the fact that they have the same dividing set to make the 
characteristic foliations agree.    This gives us a new manifold $(M',\xi')$.

The following proposition describes Legendrian surgery on the 4-manifold level.

\begin{prop}
Let $(M',\xi')$ be a contact manifold obtained by Legendrian surgery along $L$ in $(M,\xi)$, in
a 3-dimensional manner.  Then there exists a Stein cobordism from $(M,\xi)$ to $(M',\xi')$, 
obtained by attaching a 2-handle along $N(L)$. 
\end{prop}

\begin{proof}
We apply Lemma \ref{trivial-cobord} below to obtain a Stein cobordism $X=M\times[0,1].$ Then 
Legendrian surgery corresponds to attaching a 2-handle along $N(L)\subset M\times\{1\}$ in a Stein (resp.\ 
symplectic) manner, which yields a Stein (resp.\ symplectic) cobordism from $(M,\xi)$ to 
$(M',\xi')$.  (See Eliashberg \cite{Eliashberg90}.)
\end{proof}

\begin{lem} \label{trivial-cobord}
Let $(M,\xi)$ be a contact structure.  Then there exists a thickening  of $M$ to 
$X=M\times[0,1]$ and a Stein cobordism from $(M,\xi)$ to itself. 
\end{lem}

A proof of this fact appears in \cite{Eliashberg85}.

% ************************************************************************
\section{Open book decompositions}
% ************************************************************************

Recall an \dfn{open book decomposition} of a 3-manifold $M$ consists of a link $K,$
called the {\it binding}, 
and a fibration $f\co (M\setminus K)\to S^1$ such that each fiber $F$ in the fibration is a 
Seifert surface for $K.$ The manifold $M\setminus K$ is obtained by taking $F\times[0,1]$ with coordinates
$(x,t)$ and identifying $(x,0)\sim (\phi(x),1)$ via the monodromy map 
$\phi:F\stackrel{\sim}{\rightarrow} F$. 
Following Thurston 
and Winkelnkemper \cite{TW}, we construct a contact structure on $M$ from an open book 
decomposition: Let $\lambda$ be a primitive for an area form on $F$ and let $\lambda_t=
t\cdot \lambda+(1-t)\cdot \phi^*\lambda$, $t\in[0,1]$.  The 
1-form $\alpha=dt+\lambda_t$ is a contact 1-form on $F\times[0,1]$ which glues to give a 
contact structure on $M\setminus K$.  One easily checks that $\alpha$ extends over $K.$ If 
$(M,\xi)$ is obtained in this manner, then we say that the open book decomposition of $M$ is 
{\it adapted to} $\xi$. We now have the following recent result of Giroux \cite{Giroux}:

\begin{thm} \label{convex} 	
	Any contact structure $\xi$ on a closed 3-manifold $M$ admits an 
	open book decomposition of $M$ which is adapted to $\xi$. 
\end{thm}

The following lemma (and more importantly its converse) is due to the efforts of many people, 
beginning with the work of Loi and Piergallini \cite{LP} (also see Mori \cite{Mori} for an 
earlier effort), and recently culminating in the works of Giroux \cite{Giroux}, Akbulut-Ozbagci 
\cite{AO}, and Matveyev \cite{Mat}.

\begin{lem} If the monodromy $\phi:F\rightarrow F$ for an open book can be 
expressed as a product of positive Dehn twists, then the adapted contact structure 
is Stein fillable.
\end{lem}

\begin{proof}
If $\phi=id$, then the manifold $M_n$ is simply the connected sum 
of several copies of $S^1\times S^2$.  There is a unique tight contact structure $\xi_n$ on 
$M_n=\#_n (S^1\times S^2)$, and it is Stein fillable.  The uniqueness of $\xi_n$ on $M_n$ 
follows from the unique connect sum decomposition theorem of Colin \cite{Co97} and the 
uniqueness on $S^1\times S^2$ due to Eliashberg \cite{Eliashberg92}.

Assume $\phi$ consists of a single positive Dehn twist along $\gamma\subset F$.   Then the 
manifold $M$ is obtained from $M_n$ by a Dehn surgery along $\gamma$ with surgery coefficient 
one less than the framing induced on $\gamma$ by the fiber. But we can also make $\gamma$ a 
Legendrian curve in $F$ so that the framings given by the contact structure and the fibers 
agree.  (In other words, the twisting number of $\gamma$ relative to $F$ is zero.)  This is 
made possible by applying the Legendrian Realization Principle.  Note that 
to apply the Legendrian Realization Principle, a {\it fold} may be necessary (for 
details see \cite{H1}).  Thus $(M,\xi)$ is obtained from $(M_n,\xi_n)$ by a Legendrian surgery 
and hence is Stein fillable. Now, if $\phi$ is the product of $k>1$ positive Dehn twists, we 
perform $k$ Legendrian surgeries on different leaves, in order. 
\end{proof}

We are now ready to prove Theorem \ref{thm1}.  It should be pointed out that the strategy of 
proof is similar to the proof strategy in \cite{DG}, where it is proved that ``most'' 
universally tight contact contact structures on torus bundles over the circle are not 
(strongly) symplectically fillable.  

\begin{proof}[Proof of Theorem \ref{thm1}]
If $(M,\xi)$ is Stein fillable, then we are done by Lemma \ref{trivial-cobord}.
Therefore, let $(M,\xi)$ be a contact structure which is not Stein fillable. By Theorem 
\ref{convex}, there exists an open book decomposition for $M$ which is adapted to $\xi$. Let 
$K$ be the binding, $f\co (M\setminus K)\to S^1$ the fibering of the 
complement, $F$ the fiber, and $\phi$ the monodromy map. Since $(M,\xi)$ is not Stein 
fillable, any product decomposition of $\phi$ into Dehn twists must contain some negative Dehn 
twists. We view each Dehn twist as being done on a separate fiber. On a fiber just 
after one on which a negative Dehn twist was done along $\gamma$, we can take a parallel copy 
of $\gamma$ and perform a positive Dehn twist, which is tantamount to a Legendrian surgery.  If 
a compensatory positive Dehn twisted is added whenever there is a negative Dehn twist, then 
we will have a new monodromy map $\phi'$ with only positive Dehn twists. Of course $\phi'$ will 
define a different manifold $M'$ and a different contact structure $\xi'.$  However, since the 
difference in between the monodromy for $M$ and for $M'$ is just several positive Dehn twists, 
we can get from $(M,\xi)$ to $(M',\xi')$ by a sequence of Legendrian surgeries. 
Thus we have a Stein cobordism from $(M,\xi)$ to $(M',\xi').$  
\end{proof}

% ************************************************************************
\section{Overtwisted Contact Structures}
% ************************************************************************

In this section we prove Theorem \ref{thm2}.   The proof will be broken down into two 
propositions.

\begin{prop}   \label{prop1}
	Any overtwisted contact manifold is Stein cobordant to any overtwisted contact manifold.
\end{prop}

\begin{proof}
Let $(M_i,\xi_i),i=1,2$ be two overtwisted contact manifolds.
It is a well-known fact in 3-manifold topology that we can find a link $L$ in $M_1$ such
that a certain integer Dehn surgery on $L$ will yield $M_2.$ Thus we can construct a
topological cobordism $X$ from $M_1$ to $M_2$ by attaching 2-handles with the appropriate framing 
to $M_1\times[0,1].$ Moreover, one can adapt the
proof of Lemma~4.4 in \cite{Gompf98} to show that we may assume that $X$ has an almost
complex structure with complex tangencies $\xi_i$ on $M_i.$ We now apply the following 
theorem of Eliashberg (Theorem 1.3.4 in \cite{Eliashberg90}):

\begin{thm}[Eliashberg] Let $(X,J)$ be a compact, almost complex (real) 4-manifold  with
boundary $\bdry X=M_2-M_1$.  Assume $M_1$ is $J$-concave, $J$ is integrable near $M_1$,
and the corresponding contact structure $(M_1,\xi_1)$ is {\em overtwisted}.  If the cobordism 
$(X,J)$ from $M_1$ to $M_2$ consists of only 2-handle attachments, then there exists a 
deformation of $J$ (rel $M_1$) to an integrable complex structure $\widetilde J$ on $X$.
\end{thm}

Using this theorem, we obtain a 
Stein structure on $X$ for which the complex tangencies on $M_1$ are $\xi_1$ and on $M_2$ are 
some contact structure $\xi'$ homotopy equivalent to $\xi_2.$ Now, we are done if $\xi'$ is 
overtwisted, since overtwisted contact structures are classified by their 2-plane field 
homotopy type \cite{Eliashberg89}. But we can easily ensure that the contact structure on $M_2$ 
is overtwisted by adding some extra Lutz twists to $(M_1,\xi_1)$ that are disjoint from the 
regions where the 2-handles are attached. 
\end{proof}

\begin{prop}   \label{prop2}
	Given a tight contact manifold $(M,\xi)$, there exists an overtwisted contact structure
	$\xi'$ on $M$ in the same homotopy class as $\xi$ and which satisfies $(M,\xi')\prec
	(M,\xi)$.
\end{prop}

\begin{proof}
Given $(M,\xi)$, take a Legendrian curve $L\subset M$ and its standard neighborhood $N(L)$.  
Choose a framing as in Section \ref{sect2} so that the slope of the dividing set of $\bdry 
N(L)$ is $\infty$. Now, identify slopes $s\in\R\cup\{\infty\}$ with their respective 
``angles'', $[\theta_s]\in \R/\pi \Z$.  In order to distinguish the different amounts of 
``wrapping around'', we will choose a lift $\theta_s\in \R$ instead.   There exists an 
exhaustion of $N(L)$ by concentric $T^2$, where the angles of the dividing curves on the tori 
monotonically increase over the interval $[{\pi\over 2},\pi)$ as the $T^2$ move towards the 
core. 

Now, let $(M,\xi')$ be the overtwisted 3-manifold obtained by performing a full Lutz twist 
along $L$.  This replaces $N(L)$ by the solid torus $N$, where the angles of the dividing 
curves of an exhaustion by tori monotonically increase over the interval $[{\pi\over 2}, 
3\pi)$.   We claim that a full Lutz twist $(M,\xi) \stackrel{L}{\rightsquigarrow} (M,\xi')$ is 
the inverse process of a sequence of Legendrian surgeries along the same core.    In fact, take 
a Legendrian curve $K$ in $(M,\xi')$ in the same isotopy class as $L$, whose standard 
neighborhood $N(K)\subset N$ has an exhausting set of tori which spans the interval 
$[3\pi-{3\pi\over 4} ,3\pi)$.  After Legendrian surgery, the new $N$ ``rotates'' in the 
interval $[{\pi\over 2}, {5\pi\over 2})$.  Repeated application (total of 4 times) of 
Legendrian surgery will get us back to $(M,\xi)$.  Note, however, that the intermediate 
manifolds are not necessarily diffeomorphic to $M$.
\end{proof}

Combining Propositions \ref{prop1} and \ref{prop2}, we immediately get Theorem \ref{thm2}.

% ************************************************************************
\section{Concave Fillings}
% ************************************************************************

In this section we prove Theorem \ref{thm3}. Before we set out on the proof, we give a  
straightforward proof of this theorem for overtwisted contact structures.
\begin{lem}\label{otconcave}
	Theorem \ref{thm3} is true for any overtwisted contact structure.
\end{lem}
\begin{proof}
Given any overtwisted contact structure $(M,\xi)$, we know by Theorem \ref{thm2}
that there is a Stein cobordism $(X,\omega)$ from $(M,\xi)$ to $(S^3,\xi_{std}).$
Let $(Y,\omega')$ be any closed symplectic 4-manifold.  Use Darboux's theorem to
excise a small standard ball around a point in $Y$ and obtain a manifold $Y'$ with concave 
boundary $(S^3,\xi_{std}).$  We then obtain a concave filling of $(M,\xi)$ by gluing 
$(X,\omega)$ to $(Y',\omega'\vert_{Y'})$.  It is clear 
that there are infinitely many choices for $(Y,\omega')$ that will yield infinitely many different 
concave fillings for $(M,\xi)$. 
\end{proof}

\begin{lem}\label{Steinconcave}
	Theorem \ref{thm3} is true for any Stein fillable contact structure.
\end{lem}
\begin{proof}
Let $(M,\xi)$ be Stein filled by $(X,\omega).$ According to Corollary 3.3 in 
\cite{LiscaMatic97}, there is a symplectic embedding of $(X,\omega)$ into a compact K\"ahler 
minimal surface $S$ of general type. If we take $Y=\overline{S\setminus X}$, then 
$(Y,\omega|_Y)$ will be a concave symplectic filling of 
$(M,\xi).$

A slight modification of the above argument will produce 
infinitely many concave fillings.  Specifically, in a small standard 3-ball $(B^3,\xi_{std})\subset 
(M,\xi)$, there exist a right-handed Legendrian trefoil knot with $tb=1$ and a linking 
Legendrian unknot with $tb<0$.  If we add 
2-handles to $X$ along these Legendrian knots, we obtain a new Stein manifold 
$(X',\omega')$.  Embed $X'$ in a compact K\"ahler surface $S$ and remove $X$ to obtain a concave 
symplectic filling $(Y',\omega')$ of $(M,\xi).$  In the layer $X'\setminus X$ in $Y'$ 
there exists a symplectically embedded torus $T$ (see \cite{GompfN}).  
%This works for the standard knot in $S^3=\partial B^4$. Now add a 1-handle to sum this into
%our case.
Let $E(n)$ be the elliptic surface obtained by
taking the normal sum \cite{GompfNS} of $n\geq 1$ copies of the rational elliptic surface along 
regular fibers.  Then consider the symplectic manifold $Y_n=E(n)\#_T Y'$, obtained by 
taking the normal sum of $Y'$ along $T$ and $E(n)$ along a regular fiber.  
These concave fillings of $(M,\xi)$ are not related by blowing up and down,
since if they were then the compact manifolds
$S_n,$ obtained from $S$ by normal summing with $E_n,$ would also be so related. 
However, this is not the case, as 
$b_2^+(S_n)=b_2^+(S)+2n$ and $b_2^+$ is unchanged by blowing up and down.
\end{proof}

Theorem \ref{thm3} now follows from Lemma \ref{Steinconcave} and Theorem \ref{thm1}. 

\noindent
\s
{\em Acknowledgments.} This work owes its beginnings to some questions raised by 
Charlie Epstein.    The 
first author gratefully acknowledges the support of an NSF Post-Doctoral 
Fellowship(DMS-9705949), Stanford University and the American Institute of 
Mathematics.    The second author thanks the NSF (DMS-0072853), the American 
Institute of Mathematics, and IHES. 

% ----------------------------------------------------------------------
%\bibliography{../../refs/references}
% ----------------------------------------------------------------------

\end{document}